\documentclass[letterpaper, 12pt]{article}
\usepackage{amsmath,amssymb,amsfonts}
\usepackage{color}
\usepackage{pstricks}
\usepackage{graphicx}
\numberwithin{equation}{section}
\newtheorem{thm}{Theorem}

{ \pagestyle{plain}

\def\cwedge{\bigcirc\kern-1.07em\wedge\ }

\newcommand{\qed}{\hfill\fbox{}\par\vspace{.2cm}}

\begin{document}

\begin{center}
{\LARGE {4-dimensional Space forms as determined by the volumes of small geodesic balls}}
\end{center}

\begin{center}
{\large JeongHyeong Park}
\end{center}

\begin{center}

Sungkyunkwan University,
 Suwon, Korea

\end{center}
\begin{abstract}
\noindent Gray-Vanhecke conjectured that the volumes of small geodesic balls could determine if the manifold is a space form, and provided a proof for 
%the manifold $M$ with 
the compact 4-dimensional manifold, and some cases. In this paper, similar results for the 4-dimensional case are obtained, 
%for the 4-dimensional case, 
based upon tensor calculus and classical theorems rather than the topological characterizations in \cite{GV}.
%In this paper, with an additional condition on Euler characteristic, we prove the conjecture holds for dimension 4. %By utilizing results on the Einstein-Ka\"ehler manifold, we obtain that the conjecture holds for 4-dimensional harmonic manifolds.
\end{abstract}
\noindent {\it Mathematics Subsect Classification (2020)} : 53C21, 53B20

\section{Introduction}\label{sec1}

Let $M$ be a Riemannian manifold.
Let $p\in M$ and $B_r (p) = \{q\in M\ |\ d(q,p) \leq r\}$  be a
geodesic ball centered at $p$ of radius $r$ and $V_M(p,r)$ be the volume of $B_r (p)$.
In this study, we investigate how the volume of a small geodesic ball determines the geometry of the manifold.

One interesting avenue of inquiry is to investigate how the volume of a small geodesic ball is related to the curvature of the underlying manifold. There is a vast literature on this subject; we refer to \cite{B, GP, GPV, Gr, GV, Wil} and the references cited therein for further details. In this context, Gray and Vanhecke \cite{GV} posed the following conjecture:\\

{\bf{Conjecture.}} Let $M$ be an $n$-dimensional Riemannian manifold and suppose that 
%for all $p\in M$ and all sufficiently small $r> 0$,
$V_M(p,r)$ is the same as that of an $n$-dimensional manifold of constant sectional curvature $c$ for all $p$ and all sufficiently small $r>0$, then $M$ is also a space of constant sectional curvature $c$.\\

% They proved in the cases of $dim M \leqq 3$, conformally flat manifold, and Einstein manifold \cite{GV}.
% The investigation aims to solve this conjecture specifically for 
% $dim M = 4$.\\

By rescaling the metric, 
%%%
%\newpage
we let $\mathbb{S}^4 : = (\mathbb{S}^4, g_{+1})$ be the sphere of constant positive sectional curvature $1$, $\mathbb{H}^4 := (\mathbb{H}^4, g_{-1})$ be the hyperbolic space of constant negative sectional curvature $-1$ and $\mathbb{T}^4$ be the flat torus. Let $\chi(M)$ be a
Euler characteristic of $M$.
In this paper, we shall prove the following Theorem.

 \begin{thm}\label{thm:10}Let ($M$, $g$) be a $4$-dimensional compact Riemannian
 manifold.\\
 (1) Let $B_r (p) \subset M$, $B_r (q) \subset \mathbb{T}^4$.
If $V_M(p,r) = V_{\mathbb{T}^4}(q,r)$
for all sufficiently small $r$ and for all $p \in M$, $q \in {\mathbb{T}^4}$, and if $\chi(M) \ge
0$, then $M$
is flat.\\
(2) Let $B_r (p) \subset M$, $B_r (q) \subset \mathbb{S}^4$. 
 If $V_M(p,r) = V_{\mathbb{S}^4}(q,r)$ for all sufficiently small $r$ and for all $p \in M$, $q \in {\mathbb{S}^4}$, and if
 % ($B_r^{\mathbb{S}^4} (q)$) Vol($B_r
% (p)$) = Vol ($B_r^{\mathbb{S}^4} (q)$) 
% for all sufficiently small $r$ and all $p \in M$, $q \in {\mathbb{S}^4}$ and if $\chi(M) \ge
% \chi{(\mathbb{S}^4)}$ and $vol(M) \leq vol({\mathbb{S}^4})$, 
$\chi(M) \geq \frac{3}{4\pi^2} vol(M, g)$,
then $(M, g)$
is a space of constant sectional curvature $1$. In addition, if
$vol (M, g) \geq vol (\mathbb{S}^4)$, then $(M, g)$ is  $\mathbb{S}^4$.\\
%isometric to $(\mathbb{T}^4, flat)$.\\
(3) Let $B_r (p) \subset M$, $B_r (q) \subset \mathbb{H}^4$.
 If $V_M(p,r) = V_{\mathbb{H}^4}(q,r)$
for all sufficiently small $r$ and for all $p \in M$, $q \in {\mathbb{H}^4}$,
and if $\chi(M) \geq
\frac{3}{4\pi^2} vol_g ({M})$, 
then $(M, g)$
is a space of constant sectional curvature $-1$, so isometric to
${\mathbb{H}}^4/\Gamma$, where $\Gamma$ is a discrete cocompact torsion-free subgroup of isometries on
 ${\displaystyle \mathbb{H} ^{4}}$.
 \end{thm}\label{thm:99}

 % The author has learned that similar results had been proved by Gray-Vanhecke (See section 10 of \cite{GV}). 
 In this paper, the author provide a simple proof using tensor calculus and classical theorems rather than using topological characterizations in \cite{GV}. 
 
%%%%%%%%%%%%%%%%%%%%%%%%%%%%%%%%%%%%%%%%%%%%%%%%%%%%%%%%%%%%%%%%%%%
\section{Preliminaries}
%%%%%%%%%%%%%%%%%%%%%%%%%%%%%%%%%%%%%%%%%%%%%%%%%%%%%%%%%%%%%%%%%%
In this section, we prepare some notations. Let $M=(M,g)$ be
an $n$-dimensional Riemannian manifold and $\mathfrak{X}(M)$ the Lie
algebra of all smooth vector fields on $M$. We denote the
Levi-Civita connection, the curvature tensor, the Ricci tensor, and
the scalar curvature of $M$ by $\nabla$, $R$, $\rho$, and $\tau$,
respectively. The curvature tensor is defined by
\begin{equation*}
    R(X,Y)Z
    =\nabla_{[X,Y]}Z -[\nabla_X,\nabla_Y]Z
\end{equation*}
for $X$, $Y$, $Z\in\mathfrak{X}(M)$. The Weyl tensor is defined by
\begin{equation*}
\begin{split}
W_{abcd} =
R_{abcd}& - \frac{1}{n-2}(\rho_{ac}g_{bd}+\rho_{bd}g_{ac} -\rho_{ad}g_{bc}-\rho_{bc}g_{ad})\\
&+\frac{\tau}{(n-1)(n-2)}(g_{ac}g_{bd}-g_{ad}g_{bc}).
\end{split}
\end{equation*}
Then, by direct computation,  we obtain
\begin{equation}\label{w norm}
\begin{split}
{{\mid W \mid ^2 = \mid R \mid ^2 - \frac{4}{n-2}\mid \rho \mid ^2 +
\frac{2}{(n-1)(n-2)}\tau^2}},
\end{split}
\end{equation}
where $|R|^2= R_{ijkl}R^{ijkl}$, $|W|^2= W_{ijkl}W^{ijkl}$ and
$|\rho|^2= \rho_{ij}\rho^{kl}$.
Let $p\in M$ and
let $G_r (p) = \{q\in M\ |\ d(p,q) = r\}$  and $B_r (p) = \{q\in M\
|\ d(q,p) \leq r\}$  be a geodesic
 sphere and a geodesic ball centered at $p$ of radius $r$, respectively. Volume of $B_r (p)$ is represented by
\begin{eqnarray*}
{\rm Vol}(\ B_r (p)) = \int_{B_r (p)} d v_g = \int_0^r  {{ \int_{G_t
(p)} d\theta dt}}
\end{eqnarray*}where $dv_g = \sqrt{\det(g_{ij})(p)} dx^1\cdots dx^n$ is the Riemannian volume element
 with respect to local coordinates $\{x^1,\cdots, x^n\}$ of $M$ around $p$ and $d \theta$ denotes the volume form on $G_r (p)$ induced from $M$.
%  Let $V_M(p, r)
% $ be the volume of a small geodesic ball of a radius $r$ about a point $p\in M$.
Gray proved the the following holds for any Riemannian manifold $M$ and any $p\in M$ \cite{GV}:
\begin{equation*}
\begin{split}
V_M(p, r) = \frac{(\pi
r^2)^{\frac{n}{2}}}{(\frac{n}{2})!}&\{1-\frac{\tau}{6(n+2)}r^2 +
\frac{1}{360(n+2)(n+4)}(-3\mid R \mid ^2 + 8\mid \rho
\mid ^2 + 5\tau^2 \\&- 18\Delta\tau) r^4 + O(r^6)\}_p.
\end{split}
\end{equation*}
Let $\tilde{\rho} = \rho_{ij}-\frac{\tau}{n}g_{ij}$ be the traceless Ricci tensor.
Then 
\begin{equation}\label{rho}
\mid \tilde{\rho} \mid ^2 = \mid \rho \mid ^2 
-\frac{1}{n} \tau^2.
\end{equation}
Then, by using \eqref{w norm} and \eqref{rho}, we may express Gray's formula
\begin{equation}\label{GV expansion}
\begin{split}
V_M(p, r) = \frac{(\pi
r^2)^{\frac{n}{2}}}{(\frac{n}{2})!}&\{1-\frac{\tau}{6(n+2)}r^2 +
\frac{1}{360(n+2)(n+4)}(-3\mid W \mid ^2 +
(8 -\frac{12}{n-2})\mid \tilde{\rho} \mid ^2\\ &+
\frac{2}{n(n-1)} \tau^2 - 18\Delta\tau) r^4 + O(r^6)\}_p,
\end{split}
\end{equation}
holds for sufficiently small $r$. \\

For the proof of the Theorem \ref{thm:10}, we refer the following theorem.

%\begin{thm}[\cite{Wolf},\cite{DC} Prop 4.4]\label{sphere}

\begin{thm}[\cite{Wolf},\cite{DC}]\label{sphere}
 Let M be a complete Riemannian manifold of even dimension with constant sectional curvature K =1, Then it is isometric to the 
sphere $\mathbb{S}^n$ or the real projective space $\mathbb{RP}^n$.
\end{thm}

% The well-known Ledger’s formula
% gives

%  \begin{eqnarray}\label{2.2}
%   \left.\frac{d^2}{d r^2}\right|_{r= 0}\ \left(\frac{\Theta_p(r)}{r^{n-1}} \right) = - \frac{1}{3} {\rm Ric_p},
%  \end{eqnarray}
%  where $\rm Ric _p$ is the Ricci curvature at $p \in M$. Since harmonic manifold is Einstein, thus
% \begin{equation*}
% Ric(M) = n-1.
%\end{equation*}

% Hyperbolic space $\mathbb{H}^4$ is the Einstein space. There is a well-known unique result up to modulo the action of the diffeomorphism group.
% \begin{thm}[\cite{BCG}]\label{hyperbolic}
% Let $M$ be a smooth compact quotient of hyperbolic manifold $\mathbb{H}^4 = SO(4,1)/SO(4)$ and let
% $g_0$ be its standard metric of constant sectional curvature. Then every
% Einstein metric $g$ on $M$ is of the form $g= c \phi^*g_0$, where
% $\phi: M \rightarrow M$ is a diffeomorphism and $c > 0$ is a constant.
% \end{thm}

\section{Proof of Theorem \ref{thm:10}}

Now, we let $M=(M,g)$ be a 4-dimensional compact Riemannian manifold. Then,
by the Chern-Gauss-Bonnet formula, it is known that Euler characteristic $\chi(M)$ of $M$ is expressed by the
following integral formula\\
\begin{equation}\label{euler}
  \chi(M)=\frac{1}{32\pi^2}\int_M \{|R|^2-4|\rho|^2+\tau^2\}dv_g,
   \end{equation}
where $|R|^2$ and $|\rho|^2$ are the square norms of the curvature
tensor and the Ricci tensor, respectively. 
By using using \eqref{w norm} and \eqref{rho}, 
 we get
\begin{equation*}
\begin{split}
&\mid \tilde{\rho} \mid ^2 = \mid \rho \mid^2 
-\frac{1}{4} \tau^2 ,
\end{split}
\end{equation*}
\begin{equation}\label{R}
\begin{split}
&\mid W \mid ^2 = \mid R \mid ^2 - 2
\mid \tilde{\rho} \mid ^2 - \frac{1}{6}\tau^2.\\
\end{split}
\end{equation}

\noindent We remark that $(M, g)$ has constant curvature if and only if $\mid W \mid ^2 = 0$ and
$\mid \tilde{\rho} \mid ^2 = 0$. \\Now, 
by assumption $V_M(p,r) = V_{\mathbb{T}^4}(q,r)$,
by taking account of \eqref{GV expansion}, we have $$\tau_{M}(p) = \tau_{\mathbb{T}^4}=0.$$
% Consequently $\tau_M$ is
% constant and the terms $\Delta\tau $ and $\tau^2$ play no role.  
%Since $W = 0$ on $\mathbb{T}^4$, 
For $n=4$, by \eqref{GV expansion}, we
obtain
\begin{equation}\label{w}
\begin{split}
-3\mid W (p) \mid ^2 + 2 \mid \tilde{\rho}(p) \mid ^2
= 0.
\end{split}
\end{equation}
By using \eqref{w norm} and \eqref{rho}, the Euler characteristic \eqref{euler} can be expressed 
\begin{equation}\label{wy}
\begin{split}
\chi(M) 
% &= \frac{1}{32\pi^2}\int_{M}\mid\mid R\mid\mid ^2 -
% 4\mid\mid {\rho} \mid\mid ^2 + \tau^2 dvol\\
% &
= \frac{1}{32\pi^2}\int_{M}\mid W \mid ^2 {{- 2\mid
\tilde{\rho} \mid ^2 + \frac{1}{6}}}\tau^2 dvol.
\end{split}
\end{equation}
%{\blue {where $c_4 = \frac{1}{3}$.}}
 {{From \eqref{wy}, by taking account of the assumption $\chi(M) \ge  0$, and $\tau_M = \tau_{\mathbb{T}^4} = 0$, then we obtain 
% \begin{equation*}
% \begin{split}
% % &\int_{T^4}\mid W_{T^4}\mid ^2 - 2\mid\tilde{\rho}_{T^4}
% % \mid ^2 + {{\frac{1}{6}}}\tau^2 
% % dvol\\
% 0 \le &\int_{M}\mid W_{M}\mid ^2 - 2\mid \tilde{\rho}_{M}
% \mid ^2 + {{\frac{1}{6}}}\tau^2 
% dvol.
% \end{split}
% \end{equation*}
% %%%%%%%%%%%%%%%%%
% Since  $\tau_M = \tau_{\mathbb{T}^4} = 0$, 
%we obtain\\
% \begin{equation*}
% \begin{split}
% \int_{\mathbb{S}^4}\tau^2_{\mathbb{S}^4} dvol \ge \int_{M}\tau^2_{M} dvol,
% \end{split}
% \end{equation*}
% and hence
\begin{equation*}
\begin{split}
% &\int_{\mathbb{T}^4}\mid W_{\mathbb{T}^4}\mid ^2 - 2\mid
% \tilde{\rho}_{T^4} \mid ^2 dvol\\
%&
0\le \int_{M}\mid W_{M}\mid ^2 - 2\mid \tilde{\rho}_{M}
\mid ^2 dvol.
\end{split}
\end{equation*}
% Since $\mid W_{\mathbb{T}^4}\mid ^2 = \mid \tilde{\rho}_{T^4}
% \mid ^2 = 0$,
By using \eqref{w}, we have
\begin{equation*}
\begin{split}
0 \le \int_{M}{{(-2 + \frac{2}{3})}}\mid \tilde{\rho}_M
\mid^2 dvol.
\end{split}
\end{equation*}
This implies that $|\tilde{\rho}_M|^2 = 0$. Consequently, $|W_M |^2 = 0$ by \eqref{w}. 
%Thus, $M$  has constant sectional curvature. 
%%%%%%%%%%%%
From \eqref{R}, since $\tau_{M} = 0$, 
$(M, g)$ is flat.}

%%%%%%%%%%%%%
{{To prove the assertion 2, 
%we note that
If volumes of small balls $V_M(p,r)$ equal spherical balls $V_{\mathbb{S}^4}(q,r)$,
by \eqref{GV expansion}, we obtain $\tau_{M}(p) = \tau_{\mathbb{S}^4}= 12.$
% By taking account of \eqref{GV expansion}, we obtain $\tau_{M}(p) = \tau_{\mathbb{S}^4}.$ 
Consequently $\tau_M$ is
constant and the terms $\Delta\tau $ and $\tau^2$ play no role.  
Since $W = \tilde\rho = 0$ on $\mathbb{S}^4$, for $n=4$, by \eqref{GV expansion}, we again
obtain
\begin{equation}\label{w3}
\begin{split}
-3\mid W_M(p) \mid ^2 + 2 \mid \tilde{\rho_M}(p) \mid ^2
= 0.
\end{split}
\end{equation}
Now, from \eqref{wy},
% volume density function $\Theta_p(r)$ on $\mathbb{S}^4$ with constant sectional curvature $1$  is $\sin^{3} r$.
% %(resp. $\mathbb{H}^4$) (resp. $-1$) (resp. $\sinh^{n-1} r$)
% {\red{If volumes of small balls equal volume of spherical balls, using the wel-known Ledger's formula \begin{eqnarray}\label{2.2}
%   \left.\frac{d^2}{d r^2}\right|_{r= 0}\ \left(\frac{\Theta_p(r)}{r^{n-1}} \right) = - \frac{1}{3} {\rm Ric_p},
%  \end{eqnarray} (see \cite{B}), we get Ricci curvature $\rm Ric_p = 3$ at $p \in M$
% is 
% \begin{eqnarray}\label{2.2}
%   \left.\frac{d^2}{d r^2}\right|_{r= 0}\ \left(\frac{\sin^{3}r}{r^{n-1}} \right) = - \frac{1}{3} (3) = -\frac{1}{3} {\rm Ric_p},
%  \end{eqnarray}.
%  %Thus, $${\rm Ric_p} = 3$ 
% % where $\rm Ric _p$ is the Ricci curvature at $p \in M$.
% Since $S^4$ is Einstein, 
}
% Since $\tau_{M}(p) = \tau_{\mathbb{S}^4}$ by \eqref{GV expansion}, thus $\tau_{M}(p) = \tau_{\mathbb{S}^4}=12$. 
% If volumes of small balls equal volumes of spherical balls, then $\tau = \tau_{+1} = 12$, and we get 

\begin{equation}\label{c}
32\pi^2 \chi(M) = \int_M (-2 + \frac{2}{3})|\tilde{\rho_M}|^2 dvol + 24~ vol(M, g).
\end{equation}
So, in \eqref{c}, by taking account of the assumption $32\pi^2\chi(M) \geq 24 vol(M,g)$, 
then we get $|\tilde{{\rho}}|^2=0$. So, $|W_M |^2 = 0$ by \eqref{w3}. Thus, $(M, g)$ has constant sectional curvature 1. By Theorem 2, %then, $(M, g)$ is $\mathbb{S}^4$.
%or $(\mathbb{RP}^4 , g_{+1})$.
%Since $|\tilde{{\rho}}|^2=0$, by \eqref{c}, 
%we obtain $32\pi^2\chi(M) = 24 vol(M,g)$. 
%the fact that $32\pi^2\chi(M) = 24 vol(M,g)$. 
%In addition, if $vol (M, g) = vol (\mathbb{S}^4)$, taking 
%$\chi(S^4)=2$, $\chi(\mathbb{R}P^4)=1$ and 
%$vol(\mathbb{S}^4, g_{+1}) = \frac{8 \pi^2}{3}$ and $vol(\mathbb{RP}^4, g_{+1}) = \frac{4 \pi^2}{3}$ into account, 
%and  
then $(M, g)$ is a round sphere $\mathbb{S}^4$.
%$$32\pi^2 \chi(M) = \int_M (-2 + \frac{2}{3})|\tilde{\rho}|^2$$
%This implies $|\tilde{\rho}_M|^2 = 0$.  Thus we get that that equality holds. 
%Then $(M, g)$ has has constant sectional curvature 1.}}
%Consequently, $|W_M |^2 = 0$ as well. By using \eqref{w}, $M$  has constant sectional curvature 1.}

% If volumes of small balls equal volumes of spherical balls, then $\tau = \tau_{+1} = 12$, and you get 
% $$32\pi^2 \chi(M) = \int_M (-2 + 2/3)\rho^0]^2 + 1/24 vol(M, g).$$}
% So if 
% (*)  $32\pi^2\chi(M) \geq 1/24 vol(M,g)$, 
% then you get equality must hold and $(M, g)$ is a round sphere. As an extra remark, you could say (*) holds if $\chi(M) \geq \chi(S^4) = 2$ and $vol(M, g) \leq vol S^4(1)$ but this is a bit less general. 
%  To prove the assertion (3), 
% if volumes of small balls equal volumes of hyperbolic balls, then $\tau = \tau_{+1} = -12$, and you get again 
% $$32\pi^2 \chi(M) = \int_M (-2 + \frac{2}{3})|\tilde{\rho}|^2 dvol + 24 vol(M, g) \ge 24 vol(M, g).$$
% So if $32\pi^2\chi(M) \geq 1/24 vol(M,g)$, then you get again equality must hold and $(M, g)$ is a compact hyperbolic manifold.
}
%By using Theorem \ref{sphere},
% this completes the proof (1) of Theorem. Similarly,
%     (3) of Theorem can be proved %by using Theorem \ref{hyperbolic}, 
    %that is, every Einstein metric is unique up to constant, and $\tau_{M} = \tau_{\mathbb{H}^4}= - 1$. 
   
 To prove assertion 3,
   we note that if $(M, g_{-1})$ is compact Riemannian manifold of negative constant sectional curvature $-1$, then $(M, g_{-1})$ is isometric to ${\mathbb{H}}^4/\Gamma$, where $\Gamma$ is a discrete cocompact torsion-free subgroup of isometries on ${\displaystyle \mathbb{H}^{4}}$ \cite{K, R}. 
   % We let $(H^4, g_{-1}):= ({\mathbb{H}}^4/\Gamma, g_{-1})$, then it is compact hyperbolic manifold of sectional curvature $-1$.
    From \eqref{euler}, the Euler characteristic  $({\mathbb{H}}^4/{\Gamma}, g_{-1})$ is given by 
    %since $|R|^2 = 24$, $|\rho|^2 = 36$ and $\tau^2 = 144$,
    \begin{equation}\label{H}
   \chi({\mathbb{H}}^4/{\Gamma}) = \frac{3}{4\pi^2} vol (H^4, g_{-1} ).
   \end{equation}
  Now, if volumes of small balls equal volumes of hyperbolic balls, then $\tau = \tau_{-1} = -12$, and we obtain again \eqref{w3} and \eqref{c}.
Similarly as in the proof of assertion (2), we can prove the assertion (3). Considering assumption $32\pi^2\chi(M) \geq 24 vol(M,g)$, 
$$32\pi^2 \chi(M) = \int_M (-2 + \frac{2}{3})|\tilde{\rho_M}|^2 dvol + 24 vol(M, g) \ge 24 vol(M, g).$$
then we have that equality hold and $(M, g)$ is a compact hyperbolic manifold of sectional curvature $-1$.
This completes the proof of Theorem 1.
}
% Thus the conjecture holds true for complex projective space ${\mathbb{CP}}^2$ and complex hyperbolic space ${\mathbb{CH}}^2$, which implies that the conjecture holds for 4-dimensional harmonic manifolds.
%%%
\section{Acknowledgements}
The author would like to thank Professors Alice Chang and Michael Anderson for helpful comments and discussions. This work was supported by the National Research Foundation of Korea (NRF) grant funded by the Korea government (MSIT) (NRF-2019R1A2C1083957).

%%%%%%%%%%%%%%%%%%%%%%%%%%%%%%%%%%%%%%%%%%%%%%%%%%%%%%%%%%%%%%%%%%%%%%%%%%%%%%%%%%%%%

\noindent
    JeongHyeong Park\\
    Department of Mathematics, Sungkyunkwan University, Suwon 16419, Korea.\\
    E-mail: parkj@skku.edu\\

\end{document}